\documentclass{amsart}
\usepackage{graphicx}
\usepackage{amssymb}
\usepackage{tikz}
\usetikzlibrary{matrix}
\newtheorem{thm}{Theorem}[section]

\newtheorem{lem}[thm]{Lemma}

\theoremstyle{definition}
\newtheorem{defn}[thm]{Definition}

\theoremstyle{remark}
\newtheorem{rem}[thm]{Remark}

\numberwithin{equation}{section}

\begin{document}
\title[divergent Fourier series]{divergent Fourier series in function spaces near  $L^1[0;1]$}%
\author{Tengiz Kopaliani}
\author{Nino Samashvili}
\author{Shalva Zviadadze}
\address[Tengiz Kopaliani]{Faculty of Exact and Natural Sciences, Javakhishvili Tbilisi State University, 13, University St., Tbilisi, 0143, Georgia}%
\address[Nino Samashvili]{Faculty of Exact and Natural Sciences, Javakhishvili Tbilisi State University, 13, University St., Tbilisi, 0143, Georgia}%
\address[Shalva Zviadadze]{Faculty of Exact and Natural Sciences, Javakhishvili Tbilisi State University, 13, University St., Tbilisi, 0143, Georgia}%
\email[Tengiz Kopaliani]{tengizkopaliani@gmail.com}%
\email[Nino Samashvili]{n.samashvili@gmail.com}%
\email[Shalva Zviadadze]{sh.zviadadze@gmail.com}

\thanks{This work was supported by Shota Rustaveli National Science Foundation  of Georgia (SRNSFG) grant no.: FR17\_589}
\subjclass[2010]{42A16, 42A20}

\keywords{Fourier series, Uniformly bounded orthonormal system, Almost everywhere convergence, Variable exponent Lebesgue space	}%

\begin{abstract}
In this paper we generalize Bochkariev's theorem, which states that for any uniformly bounded orthonormal system $\Phi$, there exists a Lebesgue integrable function such that the Fourier series of it with respect to system $\Phi$ diverge on the set of positive measure. We characterize the class of variable exponent Lebesgue spaces $L^{p(\cdot)}[0;1]$, $1<p(x)<\infty$ a.e. on [0;1], such that above mentioned Bochkarev's theorem is valid.
\end{abstract}
\maketitle
\section{Introduction}
\par After that Kolmogorov \cite{kolmogorov1}, \cite{kolmogorov2}  gave the examples of the functions in $L^1[0;1]$ with almost everywhere and everywhere divergent trigonometric Fourier series, many authors try to generalize these results by giving the examples of functions with a.e. divergent trigonometric Fourier series from narrower Orlich's spaces. The strongest result in this direction belongs to Konyagin \cite{koniagin}, he obtained the same result for the space $L\varphi(L)$, whenever $\varphi$ satisfies $\varphi(t)=o(\sqrt{\ln t/\ln\ln t})$.
\par The similar problems with respect to other orthonormal systems were considered by different authors. One of them was the problem posed by Alexits (see \cite[pp. 326]{alexits1}, \cite[pp. 287]{alexits2}) and Olevskii \cite{olevskii} about an analogue of Kolmogorov's example of a divergent trigonometric Fourier series for general orthonormal systems that are uniformly bounded.
\par The answer to this question was given by Bochkarev \cite{bochkarev}. He proved that for every given uniformly bounded orthonormal system, there exists a function in $L^1[0;1]$ that Fourier series with respect to this system diverge at every point of some set of positive measure. It turned out that there cannot be, in general, a complete analogue of Kolmogorov's example for orthonormal systems. This follows from the fact that Kazaryan \cite{kazaryan} has constructed a complete orthonormal system that is uniformly bounded and for which every Fourier series converge on some set of positive measure.
\par The variable exponent Lebesgue spaces were introduced by Orlicz. These spaces attract spacial attention in recent three decades. Various results concerning of these spaces can be found in \cite{uribefiorenza} and \cite{DHHR}. It turns out that $L^1[0;1]=\cup L^{p(\cdot)}[0;1]$ where the union is taken over all measurable $p(\cdot)$ such that  a.e. $p(x)>1$. Using this fact authors of the paper \cite{egk} provide a different point of view on the problem of a.e. divergence of trigonometric Fourier series in the subspaces of $L^1[0;1]$. Indeed any function with Fourier series that is divergent a.e. must belong to some variable exponent space $L^{p(\cdot)}[0;1]$, $1<p(x)<\infty$ a.e. In the paper \cite{egk} authors have constructed a variable exponent space $L^{p(\cdot)}[0;1]$ with $1<p(x)<\infty$ a.e., which has in common with $L^\infty[0;1]$ the property that the space of continuous functions $C[0;1]$ is a closed linear subspace in it. Moreover, the Kolmogorov's function with a.e. divergent Fourier series belongs to the $L^{p'(\cdot)}[0;1]$ space, where $p'(\cdot)$ is a conjugate function of $p(\cdot)$. Some results about convergence of Fourier series for functions from this spaces can be found in \cite{k1} and \cite{k2}.
\par Our plan for this paper is to characterize the class of variable exponent Lebesgue spaces for which an analogue of Bochkarev's theorem is valid.
\begin{defn}
\label{def_p_ln}
Let $P_{\ln}$ be a set of all functions $p:[0;1]\to[1;\infty)$ such that $p(\cdot)$ is an increasing function, $p(0)=1$, $p(t)>1$, $t\in(0;1]$ and
\begin{equation}
\label{cond_in_def_p_ln}
\liminf_{t\to0+}(p(t)-1)\ln(e/t)<\infty.
\end{equation}
\end{defn}
\par Let now state the main result:
\begin{thm}
\label{thm_main}
For any uniformly bounded orthonormal system $\Phi$ and for any $p(\cdot)\in P_{\ln}$, there exists a measure preserving transformation $\omega:[0;1]\to[0;1]$, such that in the corresponding $L^{p(\omega(\cdot))}[0;1]$ space, there exists a function whose Fourier series with respect to $\Phi$ diverge at every point of some set of positive measure.
\end{thm}
\begin{rem}
\label{rem_in_the_dual_c_is_closed}
Let us note that (\ref{cond_in_def_p_ln}) implies that
$$
\limsup_{t\to0+}\frac{(p')^*(t)}{\ln(e/t)}>0,
$$
where $p'(\cdot)$ is conjugate function of $p(\cdot)$ and $(p')^*$ denotes the decreasing rearrangement of $p'(\cdot)$.
Then by \cite[Theorem 1.1]{kopzviad} we conclude that the space of continuous functions $C[0;1]$ is closed linear subspace in the space $L^{p'(\omega(\cdot))}[0;1]$ for some $\omega:[0;1]\to[0;1]$ measure preserving transofmation.
\end{rem}
\begin{rem}
\label{rem_in_the_proof_not_necessary}
In the proof of the Theorem \ref{thm_main} it is not necessary for $C[0;1]$ to be closed subspace in $L^{p'(\omega(\cdot))}[0;1]$. However we are constructing $\omega$ such that $C[0;1]$ is closed subspace in $L^{p'(\omega(\cdot))}[0;1]$.
\end{rem}
\par According to previouse reasoning we would like to state an open problem: Let $p(\cdot)$ is an exponent such that the space of continuous functions $C[0;1]$ is closed linear subspace in $L^{p(\cdot)}[0;1]$. Whether there exists a function in associate space $L^{p'(\cdot)}[0;1]$ such that its Fourier series diverge on a set of positive measure? (An analogue problem can be stated in case of general Banach function spaces).

\section{definitions and some auxiliary results}
Given a measurable function
$p:[0;1]\rightarrow[1;+\infty),\,\,L^{p(\cdot)}[0;1]$ denotes
the set of measurable functions $f$ on $[0;1]$ such that for some
$\lambda>0$
$$
\int_{[0;1]}\left(\frac{|f(x)|}{\lambda}\right)^{p(x)}dx<\infty.
$$
This set becomes a Banach function spaces when equipped with the
norm
$$
\|f\|_{p(\cdot)}=\inf\left\{\lambda>0:\,\,
\int_{[0;1]}\left(\frac{|f(x)|}{\lambda}\right)^{p(x)}dx\leq1\right\}.
$$
\par For the given $p(\cdot),$ the conjugate exponent $p'(\cdot)$ is defined pointwise $p'(x)=p(x)/
(p(x)-1)$, $x\in[0;1]$. 
\par The associate space of
$L^{p(\cdot)}[0;1]$ contains all measurable functions $f$ such that
$$
\|f\|_{(L^{p(\cdot)})'}'=\sup\left\{\int_{[0;1]}|f(x)g(x)|dx\,:\,g\in
L^{p(\cdot)}[0;1],\,\|g\|_{p(\cdot)}\leq1\right\}<\infty.
$$
Note that in this case the associate space of $L^{p(\cdot)}[0;1]$ is
equal to $L^{p'(\cdot)}[0;1],$ and $\|\cdot\|_{(L^{p(\cdot)})'}'$ and $\|\cdot\|_{p'(\cdot)}$ are equivalent norms (see \cite{uribefiorenza}, \cite{DHHR}).

The next result is a necessary and sufficient condition for the
embedding $L^{q(\cdot)}[0;1]$ $\subset L^{p(\cdot)}[0;1]$
(see \cite{uribefiorenza}, \cite{DHHR}).

\begin{lem}
\label{proposition_embeding_Lp_Lq}
Given the exponents $p(\cdot), q(\cdot)$, Then
$L^{q(\cdot)}[0;1]\subset L^{p(\cdot)}[0;1]$ if and only if
$p(t)\leq q(t)$ almost everywhere. Furthermore, in this case
we have
\begin{equation}
\label{estimate_embeding}
\|f\|_{p(\cdot)}\leq 2\|f\|_{q(\cdot)}.
\end{equation}
\end{lem}

\par Let $\Phi=\{\varphi_n\}_{n\in\mathbb{N}}$ be an orthonormal system (ONS) on [0;1]. For the function $f\in L^1[0;1]$  the series
$$
\sum_{n=1}^\infty c_n(f)\varphi_n(x),
$$
is called a Fourier series of $f$, where 
$$
c_n(f)=\int_0^1f(t)\varphi_n(t)dt, \:\: n=1,2,...,
$$
denotes the Fourier coefficients of $f$ with respect to the system $\Phi$.
We call ONS $\Phi$ uniformly bounded if 
$$
|\varphi_n(t)|<M, \quad\textnormal{for almost all}\: t\in[0;1] \:\textnormal{and}\: n\in\mathbb{N}.
$$
\par For any $N\in\mathbb{N}$ we denote by $D^N$ the subset of $[0;1]^N$ consisting of points $\theta = (\theta_1,...,\theta_N)$, $\theta_i\in(0;1)$, $i\in\{1,...,N\}$, that satisfy
\begin{equation}
\label{defn_theta_points_by_system}
\lim_{t\to0+}\frac{1}{t}\int_{\theta_i}^{\theta_i+t}\varphi_n(x)dx=\varphi_n(\theta_i),\quad i\in\{1,...N\},\quad n\in\mathbb{N}.
\end{equation}
Then $m(D^N)=m([0;1]^N)=1$, where $m(A)$ denotes the Lebesgue measure of the set $A$.

\par Our construction of function in Theorem 1.2 is based on the well known scheme of construction analogous function in space $L^1[0;1]$ (see \cite[Ch. IX, Lemma 1, pp. 322]{kashinsahakian}).

\begin{lem}
\label{lem_bochkarev_1}
Let the uniformly bounded ONS $\Phi$ has the property that for each $i\in\mathbb{N}$, there is a function $g_i(x)$ with $||g_i||_1\leq 2^{-i}$ that satisfies 
\begin{equation}
\label{result_bochkarevs_lemma_1}
m(\{x\in(0;1)\::\:S^*(g_i,x)>2^i\})>\alpha>0,
\end{equation}
where
$$
S^*(g_i,x)=\sup_{1\leq N<\infty}|S_N(g_i,x)|,
$$
in the last inequality, the constant $\alpha$ is independent of $i$; $S_N(g_i,x)$ are the partial sums of the Fourier series of $g_i$ with respect to the system $\Phi$. 
\par Then there is a sequence $i_1 < i_2 <...$, such that function
$$
f(x)=\sum_{k=1}^\infty g_{i_k}(x)
$$
belongs to the $L^1[0;1]$ and for which the Fourier series diverges unboundedly at each point of a set of positive measure.
\end{lem}

\par We shall state the fundamental lemma by Bochkarev (see \cite[Inequality (110), pp. 330]{kashinsahakian}) that we use in proving the main theorem.
\begin{lem}
\label{lem_bochkarev_2}
Let $\Phi$ be uniformly bounded ONS. For all natural number $N$, there exists $(\theta_1^{(N)},...,\theta_N^{(N)})\in D^N$ such that for some absolute positive constants $c_1$ and $c_2$ we have
\begin{equation}
\label{result_bochkarev_lemma}
m\left(\left\{t\in[0;1]\::\: \sup_{1\leq m <s_0N}\sum_{n=1}^m\varphi_n(t)\left[\frac{1}{N}\sum_{i=1}^N\varphi_n(\theta_i^{(N)})\right]>c_1\ln N\right\}\right)\geq c_2,
\end{equation}
for sufficient large $s_0=s_0(N,\Phi)$.
\end{lem}

\section{Poof of Theorem \ref{thm_main}}
For all $t\in[0;1]$ define the function $h(t)=\min\{p'(t), \ln(e/t)\}$. It is obvious that in this case holds
$$
\limsup_{t\to0+}\frac{h(t)}{\ln(e/t)}>0.
$$
Then there exists a sequence $t_k\downarrow0$, such that
\begin{equation}
\frac{h(t_k)}{\ln(e/t_k)}\geq a,\quad k\in\mathbb{N},
\end{equation}
for some positive number $a$. 
\par It is obvious that we can choose subsequence $t_{k_n}$ such that $2t_{k_{n+1}}<t_{k_n}$. Let's choose a positive number $c$ such that $c>e^{1/a}$, then we get 
\begin{equation}
\label{estim_c_pow_p_rearrangement_infty}
\int_0^1 c^{h(t)}dt>\int_{t_{k_{n+1}}}^{t_{k_n}}c^{a\cdot\ln(e/t_{k_n})}dt=
\end{equation}
$$
=(t_{k_n}-t_{k_{n+1}})\cdot e^{a\cdot\ln c\cdot\ln(e/t_{k_n})}>\frac{t_{k_n}}{2}\cdot\left(\frac{e}{t_{k_n}}\right)^{a\cdot\ln c}\to+\infty,\quad n\to+\infty.
$$
Accroding to (\ref{estim_c_pow_p_rearrangement_infty}) and the fact that $t_k\downarrow0$ we can choose the subsequence $(t_{k_{n_m}})$ from $(t_{k_n})$ such that 
$$
\int_{t_{k_{n_{m+1}}}}^{t_{k_{n_m}}}c^{h(t)}dt\geq1,\quad m\in\mathbb{N}.
$$
So without loss of generality we can assume that sequence $(t_k)$ is already such that
\begin{equation}
\label{properties_of_t_k}
1<a\ln(e/t_1), \quad 2t_{k+1}<t_k,\quad \int_{t_{k+1}}^{t_k}c^{h(t)}dt\geq1, \quad k\in\mathbb{N}.
\end{equation}

\par Let $\{l_k\::\:k\in\mathbb{N}\}$ be a fixed dense set on $(0;1)$ (below we will choose $l_k$ by using system $\Phi$). The density of the set $\{l_k\::\:k\in\mathbb{N}\}$ is not necessary to proof this theorem. We use this condition to construct the space $L^{p'(\omega(\cdot))}[0;1]$ such that the space $C[0;1]$ is closed subspace in it (see Remark \ref{rem_in_the_proof_not_necessary}). Let $r_k$, $k\in\mathbb{N}$ is the following numeration of the table
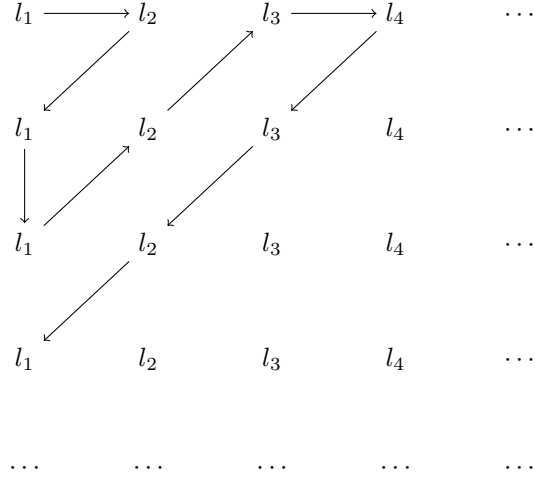
\begin{center}
\begin{tikzpicture}
\matrix(m)[matrix of math nodes,column sep=1cm,row sep=1cm]{
	l_1 & l_2 & l_3 & l_4 & \cdots \\
	l_1 & l_2 & l_3 & l_4 & \cdots \\
	l_1 & l_2 & l_3 & l_4 & \cdots \\
	l_1 & l_2 & l_3 & l_4 & \cdots \\
	\cdots & \cdots & \cdots & \cdots & \cdots\\
};
\draw[->]
(m-1-1)edge(m-1-2)
(m-1-2)edge(m-2-1)
(m-2-1)edge(m-3-1)
(m-3-1)edge(m-2-2)
(m-2-2)edge(m-1-3)
(m-1-3)edge(m-1-4)
(m-1-4)edge(m-2-3)
(m-2-3)edge(m-3-2)
(m-3-2)edge(m-4-1);
\end{tikzpicture}
\end{center}
It is obvious that for the each $l_k$ there exists sequence $(r_{k_m})$ such that $l_k=r_{k_m}$,  $m\in\mathbb{N}$.
Now let $\Delta_k=[t_{k+1};t_k]$, where $t_k$ are points possessing the property (\ref{properties_of_t_k}). Define $d_k=-t_{k+1}+r_k$ and $E_k=\Delta_k+d_k=[r_k;r_k+t_k-t_{k+1}]$. Let $g_k(t)=h(t) \cdot \chi_{\Delta_k}(t)$, $k\in\mathbb{N}$. Let introduce the functions $q_k(t)$ by the induction:
$$
q_1(t)=g_1(t-b_1)\chi_{[0;1]}(t),
$$
$$
q_k(t)=\left[q_{k-1}(t)(1-\chi_{\Delta_k}(t-d_k))+g_k(t-d_k)\right]\cdot\chi_{[0;1]}(t),\quad k>1.
$$

It is clear that $h(t)$ is decreasing and therefore $q_k(t)\leq q_{k+1}(t)$, for all $t\in[0;1]$ and all $k\in\mathbb{N}$. Also for all $k\in\mathbb{N}$ we have
\begin{equation}
\label{estim_int_q_k}
\int_0^1q_k(t)dt\leq\int_0^1h(t)dt\leq\int_0^1\ln(e/t)dt = 2.
\end{equation}
Now define a function
$$
\hat{q}(t)=\lim_{k\to+\infty}q_k(t),\quad t\in[0;1].
$$
It is clear that 
\begin{equation}
\label{inequality_q_q_k}
\hat{q}(t)\geq q_k(t)\geq a\ln(e/t_k), \:\:\: t\in E_k,\:k\in\mathbb{N}.
\end{equation}
By (\ref{estim_int_q_k}) we get that the function $\hat{q}(\cdot)$ is a.e. finite.
According to the construction it is clear that $\hat{q}^*(t)\leq h(t)\leq p'(t)$. It follows from the well known result (see \cite[Theorem 7.5]{BS}) that there exists measure preserving transformation $\omega:[0;1]\to[0;1]$ such that $\hat{q}(t)=\hat{q}^*(\omega(t))$. Now define $\bar{q}(\cdot)$ by $\bar{q}(t)=p'(\omega(t))$. Since $\hat{q}^*(t)\leq p'(t)$ it is obvious that $\hat{q}^*(\omega(t))\leq p'(\omega(t))$, then for all $t\in(0;1)$ we get following inequality
\begin{equation}
\label{estim_bar_p_by_tilde_p}
\hat{q}(t)\leq \bar{q}(t).
\end{equation}
\par Note that the space of continuous functions $C[0;1]$ is closed subset in the space $L^{\bar{q}(\cdot)}[0;1]$. Indeed, given $I=(a;b)\subset(0;1)$ let estimate $||\chi_I||_{\bar{q}(\cdot)}$. Since $c>1$, it is obvious that by (\ref{estim_bar_p_by_tilde_p}) the following inequality holds
$$
\int_I c^{\bar{q}(t)}dt\geq\int_I c^{\hat{q}(t)}dt.
$$
By the construction of $\hat{q}(\cdot)$ there exists number $k_0$ such that $E_{k_0}\subset I$. We have
$$
\int_I c^{\hat{q}(t)}dt\geq\int_{E_{k_0}}c^{\hat{q}(t)}dt\geq\int_{E_{k_0}}c^{q_{k_0}(t)}dt=\int_{E_{k_0}}c^{g_{k_0}(t-d_{k_0})}dt=
$$
$$
=\int_{r_{k_0}}^{r_{k_0}+t_{k_0}-t_{k_0+1}}c^{h(t-d_{k_0})}dt=
$$
$$
=\int_{t_{k_0+1}}^{t_{k_0}}c^{h(t)}dt\geq1.
$$
By the definition of the norm in variable Lebesgue space and by the above estimations we get that for all intervals $(a;b)$ we have $||\chi_{(a;b)}||_{\bar{q}(\cdot)}\geq 1/c$. By the \cite[Theorem 3.1]{egk} we get that the space of continuous function is closed subspace in $L^{\bar{q}(\cdot)}[0;1]$.

\par Consider the function $\bar{p}(\cdot)$ which is conjugate of $\bar{q}(\cdot)$. It is clear that $\bar{p}(t)=p(\omega(t))$, $t\in[0;1]$. Let $C>\ln(e/t_1)\cdot(a\ln(e/t_1)-1)^{-1}$, then by (\ref{inequality_q_q_k}) and (\ref{estim_bar_p_by_tilde_p}) it is obvious that
\begin{equation}
1<\bar{p}(t)\leq 1+\frac{C}{\ln(e/t_k)},\quad t\in E_k.
\end{equation}
By the last estimation we obtain
\begin{eqnarray}
\label{estim_chi_e_k_norm}
||\chi_{E_k}||_{\bar{p}(\cdot)}\asymp m(E_k)=t_k-t_{k+1}\asymp t_k.
\end{eqnarray}
Indeed, using (\ref{estimate_embeding}) we obtain
$$
\frac{1}{2}||\chi_{E_k}||_1\leq||\chi_{E_k}||_{\bar{p}(\cdot)}\leq 2||\chi_{E_k}||_{1+\frac{C}{\ln(e/t_k)}}\asymp t_k.
$$
\par Finally, by (\ref{estimate_embeding}) and (\ref{estim_chi_e_k_norm}) we obtain
\begin{equation}
\label{estim_sum_chi_e_k_norm}
\left\|\sum_{k=1}^\infty a_k \chi_{E_k}\right\|_{\bar{p}(\cdot)}\asymp\left\|\sum_{k=1}^\infty a_k \chi_{E_k}\right\|_1.
\end{equation}
\par Recall that for each $k$ there exists a sequence of natural numbers $(k_m)$, $m\in\mathbb{N}$ such that $l_k=r_{k_m}$,  $m\in\mathbb{N}$. Thus, we can rewrite (\ref{estim_sum_chi_e_k_norm}) in the following form
\begin{equation}
\label{estim_sum_sum_chi_e_k_norm}
\left\|\sum_{k=1}^\infty\sum_{m=1}^\infty a_{k_m}\chi_{E_{k_m}}\right\|_{\bar{p}(\cdot)} \asymp 
\left\|\sum_{k=1}^\infty\sum_{m=1}^\infty a_{k_m}\chi_{E_{k_m}}\right\|_{1}.
\end{equation}

\par Let select $\{l_k\::\:k\in\mathbb{N}\}$ set by using system $\Phi$. Using Lemma \ref{lem_bochkarev_2} for all $N\in\mathbb{N}$ choose $(\theta_1^{(N)},...,\theta_N^{(N)})\in D^N$ such that (\ref{result_bochkarev_lemma}) is valid. Consider the following set $\Theta':=\{\theta_i^{(N)}\::\:i\in\{1,...,N\},\:N\in\mathbb{N}\}$. If this set is not dense in [0;1], we examine a countable set $\Theta$ such that $\Theta'\subset\Theta$ and $\Theta$ is dense in [0;1] (let us note that density of the $\Theta$ in [0;1] is not essential in our construction, we only state this condition for reason that the space $C[0;1]$ will be closed subspace in $L^{\bar{q}}[0;1]$). Let $l_k$, $k\in\mathbb{N}$ is some numeration of $\Theta$.

Let $h^N=(h^{(N)}_1,...,h^{(N)}_N)$, $h^{(N)}_i>0$, $i\in\{1,..,N\}$. Consider the function defined on $(0;1)$
\begin{equation}
\label{defn_f_function}
f_N(x)=\frac{1}{N}\sum_{i=1}^N\frac{1}{h^{(N)}_i}\chi_{(\theta_i^{(N)};\theta_i^{(N)}+h^{(N)}_i)}(x), \quad ||f_N||_{L^1}=1.
\end{equation}
By (\ref{defn_theta_points_by_system}) the Fourier coefficients of $f_N$ satisfy
$$
\lim_{h^N\to0}c_n(f_N)=\frac{1}{N}\sum_{i=1}^N\varphi_n(\theta_i^{(N)}),\quad n\in\mathbb{N},
$$
therefore
\begin{equation}
\label{estim_deviation_from_partial_sums}
\lim_{h^N\to0}\left\|S_m(f_N,x)-\sum_{n=1}^m\varphi_n(x)\left(\frac{1}{N}\sum_{i=1}^N\varphi(\theta_i^{(N)})\right)\right\|_2=0.
\end{equation}
It follows from (\ref{result_bochkarev_lemma}) and (\ref{estim_deviation_from_partial_sums}) that for all $N\in\mathbb{N}$ and $\max\{h^{(N)}_i\}<h_0=h_0(N,\Phi)$
\begin{equation}
\label{estim_measure_from_lemma}
m\left(\left\{x\in(0;1)\::\: \sum_{1\leq m\leq s_0N}S_m(f_N,x)>\frac{c_1}{2}\ln N\right\}\right)>\frac{c_2}{2}.
\end{equation}
\par For all fixed $N\in\mathbb{N}$ and fixed $\theta^{(N)}_i$, $i\in\{1,...,N\}$ there exists the sequence $r_{i_k}$, $k\in\mathbb{N}$ such that $\theta^{(N)}_i=r_{i_k}$, $k\in\mathbb{N}$. Then in the definition (\ref{defn_f_function}) we assume that $h^{(N)}_i=t_{i_k}-t_{i_k+1}$ by such way that $t_{i_k}-t_{i_k+1}<h_0$.

It remains only to apply Lemma \ref{lem_bochkarev_1}: by (\ref{estim_measure_from_lemma}), for each $i \in\mathbb{N}$, we can find $N = N(i)$ and $h = h(i) > 0$ such that $g_i(x) = (\ln N)^{-1/2}f_N(x)$ has norm $||g_i||_1<2^{-i}$ and satisfies (\ref{result_bochkarevs_lemma_1}). Using Lemma \ref{lem_bochkarev_1} and (\ref{estim_sum_sum_chi_e_k_norm}) we complete the proof.


\end{document}